\documentclass[12pt]{amsart}

\newtheorem{theorem}{Theorem}
\newtheorem{lemma}{Lemma}
\newtheorem{proposition}{Proposition}
\newtheorem{corollary}{Corollary}

\theoremstyle{definition}

\theoremstyle{remark}

\newcommand{\ltwo}{L^2({\mathbb R})}
\newcommand{\ltwod}{L^2({\mathbb R}^d)}

\newcommand{\Cal}{\mathcal}
\newcommand{\Bb}{\mathbb}

\begin{document}


\title{Geometric Aspects of Frame Representations of Abelian Groups}
\author{Akram Aldroubi}
\address{Department of Mathematics, Vanderbilt University, Nashville, TN 37240}
\email{aldroubi@math.vanderbilt.edu}
\author{David Larson}
\address{Department of Mathematics, Texas A\&M University, College Station, TX 77843-3368}
\email{larson@math.tamu.edu}
\author{Wai-Shing Tang}
\address{Department of Mathematics, National University of Singapore, 10 Kent Ridge Crescent, 119260, Republic of Singapore}
\email{mattws@nus.edu.sg}
\author{Eric Weber}
\address{Department of Mathematics, Texas A\&M University, College Station, TX 77843-3368} 
\curraddr{Department of Mathematics, University of Wyoming, Laramie, WY 82071-3036}
\email{esw@uwyo.edu}
\thanks{Submitted to Trans. Amer. Math. Soc.}
\begin{abstract}
We consider frames arising from the action of a unitary representation of a discrete countable abelian group.  We show that the range of the analysis operator can be determined by computing which characters appear in the representation.  This allows one to compare the ranges of two such frames, which is useful for determining similarity and also for multiplexing schemes.  Our results then partially extend to Bessel sequences arising from the action of the group.  We apply the results to sampling on bandlimited functions and to wavelet and Weyl-Heisenberg frames.  This yields a sufficient condition for two sampling transforms to have orthogonal ranges, and two analysis operators for wavelet and Weyl-Heisenberg frames to have orthogonal ranges.  The sufficient condition is easy to compute in terms of the periodization of the Fourier transform of the frame generators. 
\end{abstract}
\date{March 4, 2002}
\maketitle



\section{Introduction}

Frames in a separable Hilbert space provide redundant but stable expansions for the Hilbert space.  Frames arise naturally in many settings, such as sampling theory, time-scale (wavelet) analysis, and time-frequency (Weyl Heisenberg or Gabor) analysis.  The ``power of redundancy'' offered by frames has been successfully demonstrated in such settings, such as Ron and Shen constructing compactly supported wavelet frames with high approximation order (\cite{RS97b} and references), Benedetto et. al. denoising signals and other data (\cite{BL98a, BT01a}).  The focus of the present paper is motivated by another aspect of the power of redundancy which allows multiplexing of signals.

A frame for a separable Hilbert space $H$ is a sequence $X = \{x_j\}_{j \in \Bb{J}}$ such that there exist constants $0 < C_1$ and $C_2 < \infty$ such that for all $x \in H$,
\[ C_1 \|x\|^2 \leq \sum_{j \in \Bb{J}} |\langle x , x_j \rangle|^2 \leq C_2\|x\|^2.\]
There exists a (possibly non-unique) dual frame sequence $\tilde{x}_j$ such that the reconstruction formula holds:
\[ x = \sum_{j \in \Bb{J}} \langle x, x_j \rangle \tilde{x}_j = \sum_{j \in \Bb{J}} \langle x, \tilde{x}_j \rangle x_j. \]

The \emph{analysis operator} for $X$ is
\[ \Theta_{X}:H \to l^2(\Bb{J}): x \mapsto (\langle x, x_j \rangle).\]
Often one wishes to know the range of the analysis operator; for instance two frames $\{x_j\}_{j \in \Bb{J}}\subset H$ and $Y = \{y_j\}_{j \in\Bb{J}} \subset K$ are similar if and only if their analysis operators have the same range in $l^2(\Bb{J}).$  If the range of $\Theta_X$ is orthogonal to the range of $\Theta_Y$, then for all $x \in H$ and $y \in K$,
\begin{align*} 
x &= \sum_{j \in \Bb{J}} (\langle x , x_j \rangle + \langle y, y_j \rangle) \tilde{x}_j \\
y &= \sum_{j \in \Bb{J}} (\langle x , x_j \rangle + \langle y, y_j \rangle) \tilde{y}_j,
\end{align*}
provided $\tilde{x}_j$ and $\tilde{y}_j$ are the \emph{standard duals} of $x_j$ and $y_j$, respectively.  Such frames are called \emph{strongly disjoint} (see \cite{HL00a}) and this procedure is what engineers call \emph{multiplexing}.

The purpose of this paper is to determine the range of the analysis operator for a frame or a Bessel sequence that arises from the action of a unitary representation of a discrete countable abelian group.  Such frames occur in the setting of regular sampling on bandlimited functions.  In section 2, we show that the range of the analysis operator can be determined by computing which characters appear in the representation.  A Bessel sequence $\{x_j\}_{j \in \Bb{J}} \subset H$ is such that a constant $C_2$ as above exists, but there is not necessarily a constant $C_1$.  A Bessel sequence also defines an analysis operator just as a frame sequence does.  Bessel sequences occur naturally when considering a subsequence of a frame sequence, which we shall do in our investigation of wavelet and Weyl-Heisenberg frames.

The problem of sampling is well known and sampling theory in various settings is well established (see \cite{BF01a, HS99a} for overviews).  In this paper, we will restrict ourselves to sampling on totally translation invariant subspaces $V_E$, i.e. all functions $f \in L^2(\Bb{R}^d)$ such that the support of $\hat{f}$ is contained in the band (i.e. measurable set) $E \subset \Bb{R}^d$.  We shall assume that the band has finite measure.  A $d \times d$ invertible matrix $A$ is a sampling matrix for the set $E$ if the lattice $\{A\Bb{Z}^d\}$ is a \emph{set of sampling} for the space $V_E$.  By a set of sampling we mean that the norm of the function is preserved by the sampling transform $\Theta_{A}$, i.e. there are constants $C_1$ and $C_2$ such that for all $f \in V_E$,
\[ C_1\|f\|^2 \leq \|\Theta_A(f)\|^2 \leq C_2 \|f\|^2, \]
where $\Theta_A$ is defined as
\[ \Theta_A: V_E \to l^2(\Bb{Z}^d): f \mapsto (f(Az))_z. \]

The sampling transform defines a unitary representation of $\Bb{Z}^d$ on $V_E$.  We apply the results of section 2 to obtain the following result regarding sampling, which is proven in section 3:
\begin{theorem} \label{T:theorem1}
Let $A$ be a sampling matrix for the set $E$ and $B$ be a sampling matrix for the set $F$.  Define the sums
\[ m_A(\xi) := \sum_{k \in \Bb{Z}}\chi_E(A^{* -1}(\xi + k)), \ 
   m_B(\xi) := \sum_{k \in \Bb{Z}}\chi_F(B^{* -1}(\xi + k)), \ \xi \in \Bb{R}^d\]
and let $X$ and $Y$ denote the support sets of $m_A$ and $m_B$, respectively.
Then we have the following:
\begin{enumerate}
\item $\Theta_A(V_E) = \Theta_B(V_F)$, if and only if $\lambda(X \Delta Y) = 0$;
\item $\Theta_A(V_E) \perp \Theta_B(V_F)$ if and only if $\lambda(X \cap Y) = 0$;
\item $\Theta_A(V_E) \cap \Theta_B(V_F) \neq \{0\}$ if and only if $\lambda(X \cap Y) \neq 0$;
\item $\Theta_A(V_E) \subset \Theta_B(V_F)$ if and only if $X \subset Y$ modulo null sets;
\item the two conditions that a) $\Theta_A(V_E) \cap \Theta_B(V_F) = \{0\}$ and b) the angle between $\Theta_A(V_E)$ and $\Theta_B(V_F)$ is strictly positive imply that $\Theta_A(V_E)$ and $\Theta_B(V_F)$ are orthogonal.
\item the samples for $A$ and $B$ commute, i.e. the projections $P_A$ and $P_B$ onto $\Theta_A(V_E)$ and $\Theta_B(V_F)$, respectively, commute.
\end{enumerate}
\end{theorem}

In section 4, we turn to frames which arise from the ordered product of two unitary groups, i.e. $\{G_1 G_2 \Psi\}$, where $G_1$ and $G_2$ are unitary groups, and $\Psi \subset H$ is a finite set.  If $\{G_1 G_2 \Psi\}$ is a frame, then the sequences $\{G_1 \psi_j\}$ and $\{G_2 \psi_j\}$ are Bessel, to which we apply results from section 2.

We define for any $\alpha \in \Bb{R}^d$ the translation operator $T_{\alpha}$ by $T_{\alpha}f = f(\cdot - \alpha)$ and the modulation operator $E_{\alpha}f = e^{2\pi i \langle \alpha, \cdot \rangle}f(\cdot )$.  For an \emph{expansive} matrix $A$ (all its eigenvalues have modulus greater than one), we define the dilation operator $D_{A}f(\cdot ) = \sqrt{|det A|} f(A \cdot )$.

A wavelet frame $\Psi=\{\psi_1,\dots,\psi_n\}$ for $\ltwod$ has the form $\{D_{A^m}T_{Xz}\Psi: m \in \Bb{Z}, z \in \Bb{Z}^d \}$, where $A$ is expansive and $X$ is nonsingular.  The analysis operator for such a frame is
\[ \Theta_{\Psi}: \ltwod \to \oplus_{j=1}^{n} l^2(\Bb{Z} \times \Bb{Z}^d): x \mapsto (\langle x, D_{A^m}T_{Xz} \psi_j \rangle). \]
Likewise, a Weyl-Heisenberg frame has the form $\{E_{Al}T_{Xz}f_j: l,z \in \Bb{Z}^d, j=1,\dots,n\}$, where $F=\{f_1,\dots,f_n\}$, with analysis operator
\[ \Theta_{F}: \ltwod \to \oplus_{j=1}^{n} l^2(\Bb{Z}^d \times \Bb{Z}^d): x \mapsto (\langle x, E_{Al}T_{Xz} f_j \rangle). \]
Here $A$ and $X$ are both nonsingular (but not necessarily expansive).

By considering the Bessel sequences $\{T_{Xz}\psi_j\}$, $\{T_{Xz}f_j\}$, and $\{E_{Al}f_j\}$, we get the following results.

\begin{theorem} \label{T:W1}
Suppose $\{D_{A^m}T_{Xz} \psi_j:j=1,\dots,n\}$ and $\{D_{B^m}T_{Yz} \phi_j:j=1,\dots,n\}$ are affine frames.  Define the sums
\[ m_j(\xi) := \sum_{k \in \Bb{Z}^d} |\hat{\psi}_j(X^{*-1}(\xi + k))|^2 \text{ and } n_j(\xi) := \sum_{k \in \Bb{Z}^d} |\hat{\phi}_j(Y^{*-1}(\xi + k))|^2. \]
Let $E_j$ and $F_j$ denote the support sets of $m_j$ and $n_j$, respectively.  The following statements hold:
\begin{enumerate}
\item $\Theta_{\Psi}(\ltwod) \subset \Theta_{\Phi}(\ltwod)$ only if $\lambda(E_j \setminus F_j) = 0$ for $j=1,\dots,n$; 
\item $\Theta_{\Psi}(\ltwod) = \Theta_{\Phi}(\ltwod)$ only if $\lambda(E_j \Delta F_j) = 0$ for $j=1,\dots,n$;
\item $\Theta_{\Psi}(\ltwod) \perp \Theta_{\Phi}(\ltwod)$, i.e. $\{D_{A^m}T_{Xz} \psi_j:j=1,\dots,n\}$ and $\{D_{B^m}T_{Yz} \phi_j:j=1,\dots,n\}$ are strongly disjoint, if $\lambda(E_j \cap F_j) = 0$ for $j=1,\dots,n$.
\end{enumerate}
\end{theorem}

\begin{theorem} \label{T:WH1}
Suppose that $\{E_{Al}T_{Xz}f_j:j=1,\dots,n\}$ and $\{E_{Bl}T_{Yz}g_j:j=1,\dots,n\}$ are Weyl-Heisenberg frames.  Define the sums
\begin{align*}
m_j(\xi) := \sum_{k \in \Bb{Z}^d} |\hat{f}_j(X^{*-1}(\xi + k))|^2 &\text{ and } n_j(\xi) := \sum_{k \in \Bb{Z}^d} |\hat{g}_j(Y^{*-1}(\xi + k))|^2; \\
\tilde{m}_j(\xi) := \sum_{k \in \Bb{Z}^d} |f_j(A^{*-1}(\xi + k))|^2 &\text{ and } \tilde{n}_j(\xi) := \sum_{k \in \Bb{Z}^d} |g_j(B^{*-1}(\xi + k))|^2.
\end{align*}
Let $E_j$, $F_j$, $\tilde{E}_j$, and $\tilde{F}_j$ denote the support sets of $m_j$, $n_j$, $\tilde{m}_j$, and $\tilde{n}_j$, respectively.  Then the following hold:
\begin{enumerate}
\item $\Theta_{F}(\ltwod) \subset \Theta_{G}(\ltwod)$ only if $\lambda(E_j \setminus F_j) = 0$ and $\lambda(\tilde{E}_j \setminus \tilde{F}_j) = 0$ for $j=1,\dots,n$;
\item $\Theta_{F}(\ltwod) = \Theta_{G}(\ltwod)$ only if $\lambda(E_j \Delta F_j) = 0$ and $\lambda(\tilde{E}_j \Delta \tilde{F}_j) = 0$ for $j=1,\dots,n$;
\item $\Theta_{F}(\ltwod) \perp \Theta_{G}(\ltwod)$ if $\lambda(E_j \cap F_j) = 0$ for $j=1,\dots,n$;
\item if $X^{*}A = Y^{*}B$, $\Theta_{F}(\ltwod) \perp \Theta_{G}(\ltwod)$ if $\lambda(\tilde{E}_j \cap \tilde{F}_j) = 0$ for $j=1,\dots,n$.
\end{enumerate}
\end{theorem}

Several remarks regarding the remainder of the paper.  All Hilbert spaces will be separable.  All groups will be discrete, countable, and Abelian.  We will use $\lambda$ to denote Haar measure on $G$, on $\hat{G}$ (the dual group), and consequently also for Legesgue measure on $\Bb{R}^d$ (this should cause no confusion).  We utilize the following definition of the Fourier transform: for $f \in L^1(\Bb{R}^d)$,
\[ \hat{f}(\xi) = \int_{\Bb{R}^d} e^{-2\pi i \langle \xi, x \rangle} f(x) dx \]
Therefore, the inversion formula is given by
\[ f(x) = \int_{\Bb{R}^d} e^{2\pi i \langle x, \xi \rangle} \hat{f}(\xi) d\xi \]
when the integral is valid.

\section{Representation Theory}

Let $G$ be a discrete, countable abelian group, and let $\pi: G \to B(H)$ be a unitary representation of $G$ on a separable Hilbert space $H$.  Let $\widehat{G}$ denote the dual group of $G$, i.e.~the group of characters on $G$.  Let $\lambda$ denote normalized Haar measure on $\widehat{G}$.  By Stone's theorem there exists a projection valued measure $p$ on $\widehat{G}$ such that 
\[ \pi(g) = \int_{\widehat{G}} g(\xi) dp(\xi).\]
Then, by the theory of projection valued measures, there exists a probability measure $\mu$ on $\widehat{G}$, a multiplicity function $m: \widehat{G}\to\{0,1,\dots,\infty\}$ and a unitary operator
\[ U: H \to \oplus_{j=1}^{k} L^2(F_j,\mu) \hookrightarrow L^2(F_1, \mu, \Bb{C}^k), \]
where $\widehat{G} \supset F_1 \supset F_2 \supset \cdots \supset F_k$ and $m(\xi) = \# \{F_j: \xi \in F_j\}$.  By $L^2(F_1, \mu, \Bb{C}^k)$, we mean functions on $F_1$ attaining values in $\Bb{C}^k$ which are measurable and square integrable with respect to the measure $\mu$.  Without loss of generality, the sets $F_j$ may all be taken to have non-zero $\mu$ measure.  Note that $k$ above could be infinite, in which case $\Bb{C}^k$ is replaced by $l^2(G)$.

The operator $U$ intertwines the projection valued measure on $H$ and the canonical projection valued measure on $\widehat{G}$, and can be thought of as a Fourier-like transform on $H$.  Indeed, if $x \in H$, we will denote $Ux$ by $\hat{x}$.  We will call the unitary $U$ the \emph{decomposition operator}.  Moreover, we will see that for our representation $\pi_A$, the Fourier transform on $\ltwod$ is used in computing the unitary $U$.

By utilizing these theorems, we have that the representation $\pi$ is unitarily equivalent, via the decomposition operator, to the representation
\[ \sigma: G \to \Bb{U}\left(\oplus_{j=1}^{k} L^2(F_j, \mu) \right) \]
given by 
\[ \sigma(g) = M_{g(\xi)} \]
where $M_{g(\xi)}$ is the multiplication operator with symbol $g(\xi)$.

Two representations $\pi$ and $\rho$ are unitarily equivalent if there is a unitary $V:H \to K$, called an \emph{intertwining operator}, such that $\rho(g)V = V \pi(g)$.  Suppose  $\pi$ and $\rho$ are two representations of an Abelian group $G$ on the Hilbert spaces $H$ and $K$, respectively.  Then there is a measure $\mu_{\pi}$ and $\mu_{\rho}$ associated to each representation, and a multiplicity function $m_{\pi}$ and $m_{\rho}$ associated to each representation.  Two representations are unitarily equivalent if and only if the measures are equivalent and the multiplicity functions agree up to a set of measure 0 (with respect to either measure).

Note that the cyclic subspaces generated by $( \chi_{F_1}(\xi),0,\dots, 0)$ and $( 0,\chi_{F_2}(\xi),\dots,0)$ are orthogonal.  It follows that if there is a vector $v \in H$ such that the collection $\{\pi(g) v: g \in G \}$ is a frame for $H$, then the representation $\pi$ must be cyclic, whence the multiplicity function must be bounded above by 1 $a.e. \mu$ and the decomposition operator $U$ maps $H \to L^2(F_1, \mu)$.

If a group representation has a frame vector as above, we say that the representation is a \emph{frame representation}.  A frame representation of a group must be unitarily equivalent to a subrepresentation of the (left) regular representation \cite[Theorem 3.11]{HL00a} and \cite{W02a}.  This, in turn, implies that the measure $\mu$ is equivalent to the restriction of Haar measure $\lambda$ to the set $F_1$.

\begin{lemma} \label{L:lemma2}
Let $\pi$ be a representation of the Abelian group $G$ on the Hilbert space $K$.  The vector $w \in K$ is a frame vector for $\pi$ if and only if the following two conditions hold:
\begin{enumerate}
\item the decomposition operator $U$ maps $K \to L^2(F, \lambda|_F)$,
\item there exist positive constants $C_1$ and $C_2$ such that for almost every $\xi \in F$,
\[ C_1 \leq|Uw(\xi)|^2 = |\hat{w}(\xi)|^2 \leq C_2. \]
\end{enumerate}
Moreover, the lower and upper frame bounds are given by the supremum of all such $C_1$ and infimum of all such $C_2$, respectively.  In particular, $x$ generates a tight frame under the action of $\pi$ if and only if $|\hat{w}(\xi)|^2 = C_1$ $\lambda$ $a.e.\xi$.
\end{lemma}
\begin{proof}
We have established the necessity of condition 1.  We now show, by contrapositive, the necessity of the bound in condition 2.  Suppose $C > 0$ is given and suppose that $|\hat{w}(\xi)|^2 > C$ for $\xi \in F' \subset F$, $F'$ a set of positive measure.  Then $\chi_{F'}(\xi)\hat{w}(\xi) \in L^2(F, \lambda)$, and consider the following calculation:
\begin{align*}
\sum_{g \in G} |\langle \chi_{F'}(\xi), \widehat{\pi(g)} \hat{w}(\xi) \rangle|^2 &= \left| \int_{\widehat{G}} \chi_{F'}(\xi) \overline{g(\xi) \hat{w}_i(\xi)} d\lambda \right|^2 \\
&= \|\chi_{F'} \hat{w}\|^2,
\end{align*}
since $G$ forms an orthonormal basis of $L^2(\widehat{G}, \lambda)$ \cite[Corollary 4.26]{F}.  We have:
\begin{align*}
\|\chi_{F'} \hat{w}\|^2 &= \int_{\widehat{G}} |\chi_{F'}(\xi) \hat{w}(\xi)|^2 d\lambda \\
&> C \int_{\widehat{G}} |\chi_F(\xi)|^2 d\lambda = C \|\chi_F\|^2,
\end{align*}
whence upper frame bound is greater than $C$.  An analogous calculation shows that if there is a set $F''$ such that $|\hat{w}(\xi)|^2 < C$ for $\xi \in F''$, then the lower frame bound is less than $C$.

We now establish the sufficiency of conditions 1 and 2.  Note that by condition 1, the collection $\{\chi_F(\xi)g(\xi): g \in G\}$ is a normalized tight frame for $L^2(F,\lambda|_F)$, since it is the image of the orthonormal basis $\{g(\xi): g \in G\}$ under the projection $M_{\chi_F}$.  Moreover, by condition 2, for all $x \in K$, $\hat{x}(\xi)\overline{\hat{w}(\xi)} \in L^2(F,\lambda|_F)$.  Consider:
\begin{align*}
\sum_{g \in G} |\langle x, \pi(g) w \rangle|^2 &= \left| \int_{\widehat{G}} \hat{x}(\xi) \overline{g(\xi) \hat{w}(\xi)} d\lambda \right|^2 \\
&= \|\hat{x}\overline{\hat{w}} \|^2 \\
&= \int_{\widehat{G}} |\hat{x}(\xi)\overline{\hat{w}(\xi)}|^2 d\lambda(\xi) \\
&\leq C_2 \int_{\widehat{G}} |\hat{x}(\xi)|^2 d\lambda(\xi) \\
&= C_2 \|x\|^2.
\end{align*}
The lower bound follows similarly.
\end{proof}
See \cite{P01a} for a similar result.

\begin{lemma} \label{L:framevectors}
Let $\pi$ be a frame representation of $G$ on $H$, with $x$ and $y$ frame vectors for $\pi$.  Then the ranges of the analysis operators $\Theta_x$ and $\Theta_y$ are identical.
\end{lemma}
\begin{proof}
We shall show that the frames generated by $x$ and $y$ are similar.  By lemma \ref{L:lemma2}, the decomposition operator $U$ maps $H$ to $L^2(F, \lambda|_F)$; and the moduli of $\hat{x}$ and $\hat{y}$ are uniformly bounded above and below.  Define the operator $M:L^2(F, \lambda|_F) \to L^2(F, \lambda|_F)$ given by 
\[ f(\xi) \to \frac{\hat{y}(\xi)}{\hat{x}(\xi)}f(\xi). \]
Clearly, $M$ is well defined and invertible and maps $\hat{x}(\xi)$ to $\hat{y}(\xi)$.  Moreover $M$ commutes with the multiplication operators $g(\xi)$.  Therefore, the operator $U^{*} M U$ is the necessary similarity.
\end{proof}

For a frame representation of $G$, there is actually a second unitary representation of $G$ in the space of coefficients.  If $\pi$ is a frame representation of $G$ on $H$, with frame vector $v \in H$, then the analysis operator $\Theta_v : H \to l^2(G)$ intertwines the operators $\pi(g)$ and $L_g$, i.e. $\Theta_v \pi(g) = L_g \Theta_v$ where $L_g: l^2(G) \to l^2(G): x(h) \mapsto x(g^{-1} h)$.  Indeed, if $e_h$ denotes the characteristic function of $\{h\}$,
\begin{align*}
\Theta_v \pi(g) x &= \sum_{h \in G} \langle \pi(g)x, \pi(h)v \rangle e_h \\
    &= \sum_{h \in G} \langle x, \pi(g^{-1}h)v \rangle e_h \\
    &= \sum_{h \in G} \langle x, \pi(h)v \rangle e_{gh} \\
    &= L_g \sum_{h \in G} \langle x, \pi(h)v \rangle e_h \\
    &= L_g \Theta_v x.
\end{align*}

Denote by $J$ the image of the analysis operator, and let $\rho$ denote the representation of $G$ on $J$ by $L_g|_J$.  The analysis operator $\Theta_v: H \to J$ is invertible, whence by the polar decomposition of $\Theta_v$, there is a unitary operator $Q:H \to J$ which intertwines the representations.  Therefore, there is a decomposition operator $\tilde{U}:J \to L^2(F', \lambda|_{F'})$, but since this is unitarily equivalent to $\pi(g)$, $F' = F$.  We have proven the following lemma.

\begin{lemma} \label{L:lemma3}
A frame representation $\pi$ of $G$ is unitarily equivalent to the associated representation $\rho$ given by the restriction of the left regular representation to the range of the analysis operator.
\end{lemma}

\begin{lemma} \label{L:lemma4}
Let $\pi_H$ and $\pi_K$ be two frame representations of $G$ on $H$ and $K$, respectively.  The ranges of the corresponding analysis operators are identical if and only if the multiplicity functions agree modulo null sets.
\end{lemma}
\begin{proof}
If the multiplicity functions agree, then both representations are unitarily equivalent to $L^2(F,\lambda|_F)$, whence by lemmas \ref{L:lemma2} and \ref{L:framevectors} the ranges will coincide.  Conversely, if the ranges coincide, then both representations are equivalent to the same representation in the coefficient space, whence they are equivalent to each other.
\end{proof}

Two frame sequences $\{x_j\}_{j \in \Bb{J}} \subset H$ and $\{y_j\}_{j \in \Bb{J}} \subset K$ are said to be \emph{similar} if there is an invertible operator $S:H \to K$ such that $Sx_j = y_j$.  Two frame sequences are \emph{strongly disjoint} if there are two frame sequences $\{x^{*}_j\}_{j \in \Bb{J}} \subset H$ and $\{y^{*}_j\}_{j \in \Bb{J}} \subset K$ with the property that $\{x^{*}_j \oplus y^{*}_j\}$ is a frame for $H \oplus K$ and $\{x_j\}$ is similar to $\{x^{*}_j\}$ and $\{y_j\}$ is similar to $\{y^{*}_j\}$.  

It can be shown that two frames are similar if and only if the ranges of their corresponding analysis operators coincide.  On the other hand, two frames are strongly disjoint if and only if the ranges of their respective analysis operators are orthogonal \cite[Theorem 2.9]{HL00a}.  Two frame representations $\pi_H$ and $\pi_K$ of an Abelian group $G$ are said to be \emph{strongly disjoint} if the representation $\pi_H \oplus \pi_K$ is also a frame representation.  This is equivalent to the frames $\{\pi_H(g) x : g \in G \}$ and $\{ \pi_K(g) y : g \in G \}$ being strongly disjoint.

\begin{lemma} \label{L:lemma4a}
Suppose $\pi_H$ and $\pi_K$ are frame representations of the Abelian group $G$.  Then the corresponding frames are strongly disjoint if and only if the multiplicity functions have disjoint support.
\end{lemma}
\begin{proof}
By lemma \ref{L:lemma2}, there exist intertwining operators
\[ U:H \to L^2(E, \lambda|_E); \qquad V:K \to L^2(F,\lambda|_F), \]
where the corresponding multiplicity functions are
\[ m_H (\xi) = \chi_E(\xi); \qquad m_K(\xi) = \chi_F(\xi). \]
We now construct a representation of $G$ on $H \oplus K$ in the usual fashion:
\[ \pi(g) = \pi_H(g) \oplus \pi_K(g); \]
we need to check whether $\pi$ is a frame representation.
It follows from the multiplicity theory that the multiplicity function $m$ for $\pi$ is given by 
\[ m(\xi) = m_H (\xi) + m_K(\xi) = \chi_E(\xi) + \chi_F(\xi) \]
and also that the measure $\mu$ is equivalent to $\lambda|_E + \lambda|_F$.
Thus, if $\lambda(E \cap F) \neq 0$, the multiplicity function attains the value 2 on a set of non-zero measure; therefore there is no single frame vector for the direct sum space.  Hence, the original frames cannot be strongly disjoint.

On the other hand, if $\lambda(E \cap F) = 0$, then the multiplicity function is at most 1, whence the representation $\pi$ is cyclic.  Moreover, the measure $\mu$ is equivalent to $\lambda|_{E \cup F}$.  Therefore, 
\[ H \oplus K \simeq L^2(E \cup F, \lambda|_{E \cup F}); \]
so by lemma \ref{L:lemma2}, $\chi_E(\xi) + \chi_F(\xi)$ is a frame vector for $\pi$.  Since there is a frame for the direct sum space, it follows that the original frames were strongly disjoint.
\end{proof}

\begin{lemma} \label{L:lemma6}
Suppose $\pi_H$ and $\pi_K$ are frame representations.  Then the range of the analysis operator for a frame vector $v$ of $\pi_H$ contains the corresponding range of the analysis operator for a frame vector $w$ of $\pi_K$ if and only if $E \supset F$.
\end{lemma}
\begin{proof}
If $E \supset F$, then it follows that $\pi_K$ is equivalent to a subrepresentation of $\pi_H$, i.e. there is an isometry $W:K \to H$ which intertwines the representations.  It follows that $\Theta_v(H) \supset \Theta_w(K)$.

Conversely, let $P_H$ and $P_K$ be the projections onto the ranges of $\Theta_v(H)$ and $\Theta_w(K)$, respectively.  Let $\tilde{U}$ and $\tilde{V}$ be the decomposition operators, respectively.  The projections $\tilde{U} P_H \tilde{U}^{*}$ and $\tilde{V} P_K \tilde{V}^{*}$ both commute with the multiplication operators $M_{g(\xi)}$.  Therefore, they are multiplication operators by characteristic functions of sets, precisely $\chi_E$ and $\chi_F$ respectively.  Now, if $E \subset F$, then $P_H P_K = P_K P_H = P_K$, whence $E \supset F$.
\end{proof}

\begin{lemma} \label{L:lemma5}
Suppose $\pi_H$ and $\pi_K$ are frame representations.  The ranges of the analysis operators for the frames have non-trivial intersection if and only if the support of the multiplicity functions have non-trivial intersection.
\end{lemma}
\begin{proof}
If the support of the multiplicity functions have non-trivial intersection, then for any vector $x \in H$ such that $supp(Ux(\xi)) \subset E \cap F$, there exists a vector $y \in K$ such that $Vy(\xi) = Ux(\xi)$ on $E \cap F$.  It follows that $\Theta_H x = \Theta_K y$.

Conversely, if the ranges of the analysis operators have non-trivial intersection, then the representation $L_g$ restricted to the intersection is non-trivial.  As in the proof of the previous lemma, the projection onto the intersection $P_D$ is given by $\chi_D$ for some $D \subset \widehat{G}$ when conjugated by $\tilde{U}$.  It follows that $D = E \cap F$.  Therefore, if the representation on intersection is non-trivial, then $E \cap F$ must have non-zero measure.
\end{proof}

We restate the previous lemmas together in the following theorem regarding frame representations of an Abelian group $G$.
\begin{theorem} \label{T:framerep}
Suppose $\pi_H$ and $\pi_K$ are frame representations of $G$ on $H$ and $K$, respectively.  Let $E,F$ denote the supports of the multiplicity functions for $\pi_H$ and $\pi_K$, respectively.  Let $\Theta_H$ and $\Theta_K$ denote the analysis operators for some frame vectors for $H$ and $K$, respectively.  Then the following hold:
\begin{enumerate}
\item $\Theta_H(H) = \Theta_K(K)$, if and only if $\lambda(E \Delta F) = 0$;
\item $\Theta_H(H) \perp \Theta_K(K)$ if and only if $\lambda(F \cap E) = 0$;
\item $\Theta_H(H) \cap \Theta_K(K) \neq \{0\}$ if and only if $\lambda(E \cap F) \neq 0$;
\item $\Theta_H(H) \subset \Theta_K(K)$ if and only if $E \subset F$ modulo null sets.
\end{enumerate}
\end{theorem}

The following two lemmata concern cyclic subrepresentations of a multiple of the regular representation.  The proof of the first may be found in \cite{B00a} (see Proposition 2.1 and Theorem 2.2).
\begin{lemma} \label{L:subrep}
Suppose that the representation $\pi$ on $H$ is a subrepresentation of some multiple of the regular representation.  If $x \in H$, then the multiplicity function associated to the subrepresentation $\pi_K$ on the cyclic subspace $K = \overline{span}\{\pi(g)x: g \in G \}$ is given by
\[ m_K(\xi) = \chi_F(\xi) \quad \text{ where } \quad F = \{ \xi \in \widehat{G}: \hat{x}(\xi) \neq 0 \}. \]
\end{lemma}

\begin{lemma} \label{L:lemma9}
Let $\pi$ be a representation on $H$ which is equivalent to a subrepresentation of some multiple of the regular representation.  Let $x \in H$ and let $E = supp(\hat{x}(\xi))$.  Suppose that there are positive constants $C_1$ and $C_2$ such that $C_1 \leq |\hat{x}(\xi)|^2 \leq C_2$ for $\lambda$-almost every $\xi \in E$.  Then $\{\pi(g)x:g \in G\}$ is a frame for its closed linear span.
\end{lemma}
\begin{proof}
Combine lemmas \ref{L:subrep} and \ref{L:lemma2}.
\end{proof}

We now turn to a similar analysis for Bessel vectors for a frame representation.  If $\pi$ is a frame representation of $G$ on $H$, then $w \in H$ is a Bessel vector for $\pi$ if there exists a constant $C_2$ such that for all $v \in H$,
\[ \sum_{g \in G} |\langle v, \pi_g w \rangle|^2 \leq C_2\|v\|^2. \]
For a Bessel vector $w$, we again have an analysis operator $\Theta_w:H \to l^2(G)$, defined as for a frame vector, only now $\Theta_w$ may not be one to one, or even have closed range.  We have the following lemma, which is an immediate corollary of lemma~\ref{L:lemma9}.

\begin{lemma}
Suppose $\pi$ is a frame representation on $H$.  Then $w \in H$ is a Bessel vector for $\pi$ provided there exists a constant $C_2$ such that $|\hat{w}(\xi)| \leq C_2$ $\lambda \ a.e. \xi$.
\end{lemma}

\begin{lemma}
Suppose $\pi$ is a frame representation on $H$.  Suppose $w \in H$ is a Bessel vector for $\pi$, and $v \in H$ is a frame vector for $\pi$.   Then $\Theta_w(H) \subset \Theta_v(H)$.
\end{lemma}
\begin{proof}
Let $h \in H$, we shall construct an $h' \in H$ such that $\langle h, \pi_g w \rangle = \langle h', \pi_g v \rangle$.  By lemma \ref{L:framevectors}, we may choose the frame vector $v$, which will be specified below.  We have
\[ \langle h , \pi_g w \rangle = \int_{\widehat{G}} \hat{h}(\xi) \overline{g(\xi)\hat{w}(\xi)} d\xi = \int_{E} \hat{h}(\xi) \overline{g(\xi)\hat{w}(\xi)} d\xi.\]
Note that since $w$ is a Bessel vector, $\hat{w}(\xi) \in L^{\infty}(\widehat{G},\lambda)$.  We let $v = U^{*}\chi_E$ and let $h' = U^{*}\hat{h}\overline{\hat{w}}$.  Then
\[ \langle h', \pi_g v \rangle = \int_{\widehat{G}} \hat{h}(\xi) \overline{\hat{w}(\xi)} \overline{g(\xi)} \chi_E(\xi) d\xi = \int_{E} \hat{h}(\xi) \overline{\hat{w}(\xi)} \overline{g(\xi)} d\xi.\]
\end{proof}

\begin{lemma} \label{L:lemmaBessel}
Suppose $\pi,\rho$ are frame representations of $G$ on $H,k$, respectively. Suppose $v \in H$ and $w \in K$ are Bessel vectors for $\pi$ and $\rho$.  Let $E$ and $F$ denote the supports of $\hat{v}(\xi)$ and $\hat{w}(\xi)$, respectively.  The following hold:
\begin{enumerate}
\item $\Theta_{v}(H) \subset \Theta_{w}(K)$ only if $\lambda(E \setminus F) = 0$;
\item $\Theta_{v}(H) = \Theta_{w}(K)$ only if $\lambda(E \Delta F) = 0$;
\item $\Theta_{v}(H) \perp \Theta_{w}(K)$ if $\lambda(E \cap F) = 0$.
\end{enumerate}
\end{lemma}
\begin{proof}
We prove item 1 by contrapositive.  Let $E' = E \setminus F$, $\lambda(E') \neq 0$.  Consider the vector $x_{0} = U^{*}\chi_{E'} \in H$; we compute $\Theta_{v}(x_{0})$
\[
\langle x_{0}, \pi_g v \rangle = \int_{\widehat{G}} \chi_{E'}(\xi) \overline{g(\xi)\hat{v}(\xi)} d\xi. \]
Since $\{g(\xi)\}$ is an orthonormal basis, for $\Theta_{w}(y_{0}) = \Theta_{v}(x_{0})$, we must have that $\hat{y}_{0}(\xi) = \overline{\hat{v}(\xi)}$ a.e., which is not possible.

For item 3, we compute the multiplicity functions of the frame subrepresentations of $G$ on the subspaces generated by $v$ and $w$.  By lemma \ref{L:subrep}, the multiplicity functions are given by the supports of $\hat{v}(\xi)$ and $\hat{w}(\xi)$, whence the frame vectors in the subrepresentations are strongly disjoint if and only if the supports are disjoint.  Item 3 now follows from lemma \ref{L:lemma4a}.
\end{proof}

In \cite{HL00a}, frame sequences $\{x_n\} \subset H$ and $\{y_n\} \subset K$ are defined to be strongly disjoint if $\{x_n \oplus y_n\}$ is again a frame for $H\oplus K$.  It is then proven that they are strongly disjoint if and only if the analysis operators have orthogonal ranges.  However, with Bessel sequences, $\{x_n \oplus y_n\}$ will always be a Bessel sequence for $H \oplus K$.  Thus, we define Bessel sequences to be strongly disjoint if their analysis operators have orthogonal ranges.  A restatement, then, of the above lemma is that the Bessel sequences $\{\pi(G)v\}$ and $\{\rho(G)w\}$ are strongly disjoint if and only if their corresponding multiplicity functions have disjoint support.

\section{Sampling Theory.}

In this section, we set out to prove theorem 1.   We define the function $\phi_E \in V_E$ by $\hat{\phi}_E = \chi_E$.  The Fourier inversion formula is valid for $f \in V_E$, therefore
\[ f(Az) = \langle f(\cdot), \phi_E(\cdot - Az) \rangle,\]
whence our stability criterion translates into the condition that the collection $\{\phi_E(\cdot - Az): z \in \Bb{Z}^d\}$ forms a frame for $V_E$.

Since $V_E$ is totally translation invariant, for any matrix $A$, we have a representation $\pi_A$ of the integers $\Bb{Z}^d$ on the space $V_E$ given by $\pi_A(z) = T_{Az}$.  Therefore, our frame criterion above now becomes the question of when $\phi_E$ is a \emph{frame vector} for the representation $\pi_A$.

We remark here that others have used abstract harmonic analysis in investigating problems in sampling.  Kluvanek \cite{K65a} was the first to do so; Dodson and Beaty \cite{DB99a} has an excellent overview of the main ideas.  See also Behmard and Faridani \cite{F94a, BF00a} and Feichtinger and Pandey \cite{FP99a}.  Note that the techniques in those papers are predominantly analysis on locally compact Abelian groups, as opposed to our technique of analyzing group representations.

We begin by considering the representation of the integers $\Bb{Z}^d$ on all of $\ltwod$ given by $\sigma_{A}(z) = T_{Az}$ and computing the decomposition operator $U$ for this representation.  We shall let $\lambda$ denote Haar (Lebesgue) measure on $\Bb{T}^d$.  Note that the dual group to $\Bb{Z}^d$ is in fact $\Bb{T}^d$ via the association $\xi \to e^{-2 \pi i \langle \xi, z \rangle}$.

\begin{proposition} \label{P:decomp}
The decomposition operator for the representation $\sigma_{A}$ is
\[ U: \ltwod \to L^2(\Bb{T}^d, \lambda, l^2(\Bb{Z}^d)) \]
given by
\[ Uf(\xi)[k] = |\det A^{*}|^{-1/2} \hat{f}(A^{* -1}(\xi + k)). \]
\end{proposition}
\begin{proof}
Let us denote the square summable sequence $|\det A^{*}|^{-1/2} \hat{f}(A^{* -1}(\xi + k))$ by $\vec{\hat{f}}(\xi)$.  We first show that $U$ does the necessary intertwining, i.e.~converts the operator $T_{Az}$ into multiplication by $e^{-2 \pi i \langle \xi, z \rangle}$.
\begin{align*}
UT_{Az}f(\xi)[k] &= |\det A^{*}|^{-1/2} \widehat{T_{Az}f}(A^{* -1}(\xi + k)) \\
    &= |\det A^{*}|^{-1/2} e^{-2 \pi i \langle A^{* -1}(\xi + k), Az \rangle}\hat{f}(A^{* -1}(\xi + k)) \\
    &= |\det A^{*}|^{-1/2} e^{-2 \pi i \langle \xi + k, z \rangle}\hat{f}(A^{* -1}(\xi + k)) \\
    &= |\det A^{*}|^{-1/2} e^{-2 \pi i \langle \xi, z \rangle}\hat{f}(A^{* -1}(\xi + k)).
\end{align*}
Therefore,
\[ UT_{Az}f = e^{-2 \pi i \langle \xi, z \rangle} \vec{\hat{f}}(\xi). \]
Furthermore, it follows by Plancherel's formula that $U$ is unitary.
\end{proof}

We denote by $\pi_A$ the representation of $\Bb{Z}^d$ on $V_E$ given by $\pi(z) = T_{Az}$.  We wish to compute the multiplicity function of the representation $\pi_A$.  Let $\phi \in \ltwod$ and define
\[ V(\phi) := \overline{span}\{T_{Az} \phi: z \in \Bb{Z}^d \}, \]
to be the cyclic subspace generated by $\phi$.
\begin{corollary}
The multiplicity function associated to the subrepresentation of $\Bb{Z}^d$ on $V(\phi)$ is 
\[ m_{\phi}(\xi) = \chi_F(\xi)\]
where
\[ F= \{\xi \in \Bb{T}^d: \sum_{k \in \Bb{Z}^d} |\det A^{*}|^{-1} |\hat{\phi}(A^{* -1}(\xi + k))|^2 \neq 0 \}.\]
\end{corollary}
\begin{proof}
This follows from lemma \ref{L:subrep} and proposition \ref{P:decomp} above.
\end{proof}

We now wish to compute the multiplicity function of the representation $\pi_A$.  Note that by the preceding corollary, the support of $m_A(\xi)$ is precisely the support of
\[ \sum_{k \in \Bb{Z}^d} \chi_E(A^{* -1}(\xi + k)),\]
since the sum gives the multiplicity function associated to $V(\phi_E)$.  In order to compute all of the multiplicity function for $\pi_A$, we decompose the representation into orthogonal cyclic subrepresentations, that together, sum to the entire representation.  The multiplicity function will then be the sum of the corresponding multiplicity functions.

Let $Q_z = A^{*-1}([-1/2, 1/2)^d + z)$.  Let $E_z = E \cap Q_z$ and define $f_z \in V_E$ by $\hat{f}_z = \chi_{E_z}$.  The exponentials $\{e^{-2 \pi i \langle \xi, Aq \rangle} : q \in \Bb{Z}^d \}$ form an orthogonal basis for $Q_z$.  It follows that the cyclic subspace $V(f_z)$ has the property that $\widehat{V}(f_z) = L^2(E_z)$; whence
\[ \oplus_{z \in \Bb{Z}^d} V(f_z) = V_E. \]
Moreover, the multiplicity function for each $V(f_z)$ is 
\[ \chi_{E_z}(\xi) = \sum_{k \in \Bb{Z}^d} \chi_{E_z}(A^{*-1}(\xi + k)).\]  
Therefore, we have proven the following lemma.

\begin{lemma} \label{L:lemma8}
The multiplicity function for the representation $\pi_A$ on $V_E$ is
\[ m(\xi) = \sum_{k \in \Bb{Z}^d} \chi_{E_k}(\xi) = \sum_{k \in \Bb{Z}^d} \chi_E(A^{* -1}(\xi + k)). \]
\end{lemma}

Recall that in order for the matrix $A$ to be a sampling matrix for $V_E$, the representation $\pi_A$ must be cyclic on $V_E$, which in turn implies that the multiplicity function must attain only the values 0 and 1.  Thus, we have the following corollary.

\begin{corollary} \label{C:sampling}
The matrix $A$ is a sampling matrix for the set $E$ if and only if
\[ \sum_{k \in \Bb{Z}^d} \chi_E(A^{* -1}(\xi + k)) \leq 1 \ a.e. \ \xi. \]
Moreover, the sum yields the associated multiplicity function.
\end{corollary}

This corollary recaptures the known result on when $A$ is a sampling matrix.  See \cite{DS85a} for the case $d = 1$ and \cite[Theorem 1.4]{BF01a} for the general case.  Indeed, in the terminology of \cite[Theorem 1.4]{BF01a}, this says that the set $E$ must be a subset of a unit cell of the reciprocal lattice $A^{*-1}\Bb{Z}^d$.

\begin{proof}[\bf{PROOF of Theorem 1.}]  We see that parts 1 through 4 follow from theorem \ref{T:framerep} and lemma \ref{L:lemma8}.  It only remains to show parts 5 and 6 of theorem \ref{T:theorem1}.  

Recall from above that there is an equivalent representation of $G$ on the range of the sampling transform; call the ranges $J_A$ and $J_B$ and the decomposition operators $\tilde{U}_A$ and $\tilde{U}_B$, respectively.  Let $P_A$ and $P_B$ be the projections onto $J_A$ and $J_B$.  We have that
\begin{align*}
\tilde{U}_A P_A \tilde{U}_A^{*} &= \chi_E(\xi) \in L^{\infty}(\widehat{G}) \\
\tilde{U}_B P_B \tilde{U}_B^{*} &= \chi_F(\xi) \in L^{\infty}(\widehat{G}).
\end{align*}
Therefore, the subspaces either have non-trivial intersection or are orthogonal. Finally, by the same computation, the projections $P_A$ and $P_B$ commute, whence by corollary 2.15 in \cite{HL00a}, the samples commute.
\end{proof}

\subsection{Unions of Lattices}

We present here an idea of how to extend these results to unions of lattices.  Moreover, some of the lattices could be shifted (see \cite{BF00a}).  We wish to extend theorem \ref{T:theorem1} to the case of sampling on this type of set.  However, we need to choose a convention, and that is how to order the samples.  Thus, if $A_1,\dots,A_n$ are $d \times d$ matrices, we shall consider the sampling transform as given below:
\[ \Theta_A: V_E \to \oplus_{i=1}^{n}l^2(\Bb{Z}^d): f \mapsto (f(A_1 z_1),\dots,f(A_n z_n)). \]
Note that the sampling transform, and hence our statements, depend on the ordering chosen for the $A_i$'s.  Indeed, see our example below.

We first note that if $\{A_iz:z \in \Bb{Z}^d, i=1,\dots,n\}$ is a set of sampling for $E$, that the multiplicity function for each $\pi_{A_i}$ is not necessarily bounded by 1.  Therefore, the collection $\{\pi_{A_i}(z)\phi_E:z \in \Bb{Z}^d\}$ is not necessarily a frame for $V_E$.  However, that collection is a frame for its closed linear span.

\begin{lemma} \label{L:lemma10}
The collection $\{\pi_{A_i}(z)\phi_E:z \in \Bb{Z}^d\}$ is a frame for its closed linear span.
\end{lemma}
\begin{proof}
Note that if $|U\phi_E|^2 \neq 0$, then $|U\phi_E|^2 \geq 1$.  Moreover, the collection is certainly a Bessel set, since it is a subset of a frame; therefore $|U\phi_E|^2 \leq B$ for some constant $B$.  The statement then follows by lemma~\ref{L:lemma9}.
\end{proof}

Let $K_i \subset V_E$ be the cyclic subspace generated by $\phi_E$.  Let $\Theta_{A_i}$ denote the sampling transform for the lattice $A_i\Bb{Z}^d$.  Since the sampling transform $\Theta_{A_i}$ can be computed by inner products, for $f \in K_i^{\perp}$, $\Theta_{A_i} f = 0$.  As a result, the range of the sampling transform is determined by $K_i$, i.e. $\Theta_{A_i}V_E = \Theta_{A_i}K_i$.  Moreover, the range is closed, and is shift invariant.  The multiplicity function for the representation on the range is given by the support set of the sum
\[ \sum_{k \in \Bb{Z}^d} \chi_E(A_i^{* -1}(\xi + k)). \]

This sum parametrizes the range of the sampling transform $\Theta_{A_i}$ as before.  Hence, we proceed by inspecting the ranges of $\Theta_{A}$ and $\Theta_{B}$ coordinate wise.  Additionally, the range of $\Theta_{A_i}$ is not altered by shifting the lattice $A\Bb{Z}^d$.  Indeed, suppose the sampling $\Theta^{'}_{A_i}$ is done on the shifted lattice $Ak + k_0$ for some $k_0$.  Then it is easily seen that $\Theta_{A_i}T_{k_0}f = \Theta^{'}_{A_i}f$, whence our statement below remains valid for shifted lattices.  We say that $A_1,\dots,A_n$ are sampling matrices for $E$ if the set $\{A_iz:z \in \Bb{Z}^d, i=1,\dots,n\}$ is a set of sampling for $V_E$.
\begin{theorem} \label{T:main2}
Let $A_1,\dots,A_n$ be full rank sampling matrices for the set $E$ and $B_1,\dots,B_n$ be full rank sampling matrices for the set $F$.  Define the sums
\[ m_{A_i}(\xi) := \sum_{z \in \Bb{Z}}\chi_E(A_{i}^{* -1}(\xi + k)), \ 
   m_{B_i}(\xi) := \sum_{z \in \Bb{Z}}\chi_F(B_{i}^{* -1}(\xi + k)), \ \xi \in \Bb{R}^d\]
and let $X_i$ and $Y_i$ denote the support sets of $m_{A_i}$ and $m_{B_i}$, respectively.
Then we have the following:
\begin{enumerate}
\item $\Theta_A(V_E) = \Theta_B(V_F)$ only if $\lambda(X_i \Delta Y_i) = 0$ for $i=1,\dots,n$;
\item $\Theta_A(V_E) \perp \Theta_B(V_F)$ if $\lambda(X_i \cap Y_i) = 0$ for $i=1,\dots,n$;
\item $\Theta_A(V_E) \cap \Theta_B(V_F) \neq \{0\}$ only if $\lambda(X_i \cap Y_i) \neq 0$ for some $i$;
\item $\Theta_A(V_E) \subset \Theta_B(V_F)$ only if $X_i \subset Y_i$ modulo null sets, for $i=1,\dots,n$;
\end{enumerate}
\end{theorem}
\begin{proof}
Let us first demonstrate the necessity conditions in parts 1, 3, and 4.  Indeed, it is readily apparent that if the range of $\Theta_A$ coincides with (intersects, is contained in, respectively) the range of $\Theta_B$, then they must do so coordinate wise.  In other words, the range of $\Theta_{A_i}$ coincides with (intersects, is contained in, respectively) the range of $\Theta_{B_i}$ for all $i$.  The necessary conditions now follow by lemma \ref{L:lemma10} and the proof of theorem \ref{T:theorem1}.

Conversely, for the ranges of $\Theta_A$ and $\Theta_B$ to be orthogonal, it is sufficient that the ranges of $\Theta_{A_i}$ and $\Theta_{B_i}$ to be orthogonal individually.  Again, the ranges are described by lemma \ref{L:lemma10}, which in combination with lemmas \ref{L:lemma4} and \ref{L:lemma6} establish part 2.
\end{proof}

We remark that part 2 above can be extended (trivially) to the case when the number of sampling matrices for $\Theta_A$ and $\Theta_B$ are different.  One simply needs to pad the difficient sampling transform with 0's in the final coordinates.

\subsection{Examples}

Our first example arises from wavelet sets.  Let $A$ be an expansive dilation matrix.  A wavelet set $W$ (see for example \cite{DLS1, DLS2}) associated to $A$ is a measurable set such that
\[ \hat{\psi}_W = \chi_W \]
 is the Fourier transform of a wavelet.  If $\psi_W$ is an MRA wavelet, then (see \cite{Wan}),
\[ \sum_{j = 1}^{\infty} \sum_{k \in \Bb{Z}^d} \chi_W(A^{*j}(\xi + k)) = 1. \]
It follows that the sums
\[ \sum_{k \in \Bb{Z}^d} \chi_W(A^{*j}(\xi + k)) \quad \text{ and } \quad
   \sum_{k \in \Bb{Z}^d} \chi_W(A^{*l}(\xi + k)) \]
have disjoint support for distinct positive integers $j$ and $l$. Hence, the sampling transforms for $A^{-j}\Bb{Z}^d$ and $A^{-l}\Bb{Z}^d$ have orthogonal ranges.

Consider now the Shannon wavelet set $E = [-1, -\frac{1}{2}) \cup [\frac{1}{2}, 1)$.  We wish to investigate the uniform sampling sets for $V_E$.  Note that by the Beurling Density theorem, any set (irregular even) with Beurling Density larger than 2 will be a sampling set for $E$.  By corollary \ref{C:sampling} we wish to find all $A \in \Bb{R}$ such that 
\[ \sum_{k \in \Bb{Z}} \chi_E(A^{-1}(\xi + k)) \leq 1 \ a.e.\xi.\]
This is equivalent to the condition that $AE$ is $1$ translation congruent to a subset of $[0, 1)$.  Clearly, then, for $A \leq \frac{1}{2}$, this condition is satisfied.  Moreover, for $A =1$, the sum is actually identically $1$.  Therefore, $A$ is a sampling matrix for $E$ if  $A \in (0,\frac{1}{2}) \cup \{1\}$.

Note that
\[ \int_{0}^{1} \sum_{k \in \Bb{Z}} \chi_E(A^{-1}(\xi + k)) d\lambda = A\lambda (E). \]
Therefore, since $\lambda(E) = 1$, we must have that $A \leq 1$.  Now, for $A \in (\frac{1}{2}, 1)$, $AE = [-A, -\frac{A}{2}) \cup [\frac{A}{2}, A)$ has the property that $AE \cap (AE = 1) = (\frac{A}{2}, 1 - \frac{A}{2})$.  Hence, the sum above attains the value 2 on a non-null set.

We remark that our set $E$ is $1$ translation congruent to $[-\frac{1}{2},\frac{1}{2})$; such congruence preserves the orthonormality of the integer exponentials.  However, by our example, $1$ translation congruence does not preserve sampling sets.  Moreover, there is a gap between the necessary and sufficient densities, unlike with the set $[-\frac{1}{2}, \frac{1}{2})$.  Indeed, we have seen that a Beurling density of at least 1 is necessary, but a Beurling density of at least $2$ is sufficient, and the bound 2 is sharp.

Lastly, we consider two matrices together.  Our example demonstrates that the order chosen for the matrices affects the range of the transform.  Let $A_1 = \frac{1}{3}$ and $A_2 = \frac{2}{3}$.  We compute the sums, restricted to $[0,1)$:
\[
\sum_{k \in \Bb{Z}} \chi_E(\frac{3}{2}(\xi + k)) = 2 \chi_{[\frac{1}{3},\frac{2}{3})}(\xi), \ \sum_{k \in \Bb{Z}} \chi_E(3(\xi + k)) = \chi_{[\frac{1}{6},\frac{1}{3})}(\xi) + \chi_{[\frac{2}{3},\frac{5}{6})}(\xi).
\]
Therefore, the sampling transforms for $A_1$ and $A_2$ individually are orthogonal, whence by theorem \ref{T:main2}, together the sampling transform for $\{A_1, A_2\}$ is orthogonal to $\{A_2, A_1\}$.  Note that $A_1$ by itself is a sampling matrix for $E$.

\section{Affine, Quasi-Affine, and Weyl Heisenberg Frames}

We now apply our results to affine,quasi-affine, and Weyl Heisenberg frames.  Before proceeding, we remark that the analysis of wavelets using the spectral multiplicity methods was begun by Baggett, Medina, and Merrill \cite{BMM99a} in their study of generalized multiresolution analyses and MSF wavelets.  Other spectral methods in the analysis of wavelets can be found in Jorgensen, et al in \cite{BEJ00a, J01b}.
\begin{lemma} \label{L:waveL}
Suppose $G, K$ are discrete countable abelian groups.  Let $\pi^i$ be a unitary representation of $G$ on $H_i$, and let $\rho^i$ be a unitary representation of $K$ on $H_i$, for $i=1,2$.  Any frames of the form $\{\rho^1_k \pi^1_g \psi_j:k \in K, g \in G, j=1,\dots,n\}$ and $\{\rho^2_k \pi^2_g \phi_j: k \in K, g \in G, j=1,\dots,n\}$ are strongly disjoint if the Bessel sequences $\{\pi^1_g \psi_j\}$ and $\{\pi^2_g \phi_j\}$ are strongly disjoint for $j=1,\dots,n$.
\end{lemma}
\begin{proof}
Let $h_1 \in H_1$ and $h_2 \in H_2$.  We compute
\[
\sum_{k \in K} \sum_{g \in G} \sum_{j=1}^{n} \langle h_1, \rho^1_k \pi^1_g \psi_j \rangle \overline{\langle h_2, \rho^2_k \pi^2_g \phi_j \rangle} =
\sum_{k \in K} \sum_{g \in G} \sum_{j=1}^{n} \langle \rho^1_{k^{-1}} h_1,  \pi^1_g \psi_j  \rangle \overline{\langle \rho^2_{k^{-1}} h_2, \pi^2_g \phi_j \rangle} = 0.
\]
\end{proof}

Additionally, under certain conditions, the strong disjointness of the Bessel sequences $\{\rho^1_k \psi_j\}$ and $\{\rho^2_k \phi_j\}$ imply the strong disjointness of the frames $\{\rho^1_k \pi^1_g \psi_j \}$ and $\{\rho^2_k \pi^2_g \phi_j \}$.

\begin{proof}[\bf{PROOF of Theorem 2.}]
For item 1, if $\Theta_{\Psi}(\ltwod) \subset \Theta_{\Phi}(\ltwod)$, then for every $f \in \ltwod$, there exists a $g$ such that
\[ \langle f, T_{Xz} \psi_j \rangle = \langle g, T_{Yz} \phi_j \rangle \quad \forall z \in \Bb{Z}^d. \]
If $\Theta_{\psi_j}$ and $\Theta_{\phi_j}$ denote the analysis operators for the Bessel sequences $\{T_{Xz}\psi_j:z \in \Bb{Z}^d\}$ and $\{T_{Yz}\phi_j:z \in \Bb{Z}^d\}$ respectively, then $\Theta_{\psi_j}(\ltwod) \subset \Theta_{\phi_j}(\ltwod)$.  Therefore, by lemma \ref{L:lemmaBessel}, the multiplicity function for the representation on the cyclic subspace generated by $\psi_j$ is dominated by that for $\phi_j$.  The statement follows from lemma \ref{L:subrep} that the supports of $m_j(\xi)$ and $n_j(\xi)$ are the supports of the multiplicity functions.

Item 2 follows immediate from item 1.  The argument for item 3 is analogous to the argument for item 1 in combination with lemma \ref{L:waveL}.  
\end{proof}

The quasi-affine system was introduced by Ron and Shen \cite{RS97b} in order to study affine wavelets using techniques from shift invariant space theory.  For quasi-affine systems, typically one assumes that the expansive matrix $A$ preserves the integer lattice $\Bb{Z}^d$, though that is not required for our purposes.  Given an expansive matrix $A$, define the $L^1$ isometry
\[ \widetilde{D}_A: L^2(\Bb{R}^d) \to L^2(\Bb{R}^d): f(x) \mapsto |\det A| f(Ax). \]
The quasi affine system is given by
\[ \Cal{U}_{A,X}^q := \{D_{A^m} T_{Xz}: m \geq 0; z \in \Bb{Z}^d \} \cup \{T_{Xz} \widetilde{D}_{A^m} : m < 0; z \in \Bb{Z}^d\}. \]
If we denote the affine system by $\Cal{U}_{A,X} := \{D_{A^m} T_{Xz}: m \in \Bb{Z}; z \in \Bb{Z}^d\}$, then the fundamental result of Ron and Shen is that $\Cal{U}_{A,X}(\Psi)$ is a frame if and only if $\Cal{U}_{A,X}^q(\Psi)$ is a frame; moreover, the frame bounds are the same.  However, unlike with the affine system, where we needed only to compare the Bessel sequences $\{T_{Xz} \psi_j\}$ and $\{T_{Yz} \phi_j\}$, with the quasi-affine system we must compare the Bessel sequences $\{T_{Xz} \widetilde{D}_{A^r} \psi_j\}$ and $\{T_{Yz} \widetilde{D}_{B^r} \phi_j\}$ for $j=1,\dots,n$ and $r \leq 0$.  Let $\Theta_{\Psi}^q$ denote the analysis operator for the quasi-affine system with generators $\Psi = \{\psi_1,\dots,\psi_n\}$.

\begin{theorem}
Suppose $\Cal{U}^q(\Psi)$ and $\Cal{U}^q(\Phi)$ are quasi-affine frames.  Define the sums
\[ m_j^r(\xi) := \sum_{k \in \Bb{Z}^d} |\hat{\psi}_j(A^{*-r}X^{*-1}(\xi + k))|^2 \]
and
\[ n_j^r(\xi) := \sum_{k \in \Bb{Z}^d} |\hat{\phi}_j(B^{*-r}Y^{*-1}(\xi + k))|^2. \]
Let $E_j^r$ and $F_j^r$ denote the support sets of $m_j^r$ and $n_j^r$, respectively.  The following statements hold:
\begin{enumerate}
\item $\Theta^q_{\Psi}(\ltwod) \subset \Theta^q_{\Phi}(\ltwod)$ only if $\lambda(E_j^r \setminus F_j^r) = 0$ for $j=1,\dots,n$ and for $r \leq 0$; 
\item $\Theta^q_{\Psi}(\ltwod) = \Theta^q_{\Phi}(\ltwod)$ only if $\lambda(E_j^r \Delta F_j^r) = 0$ for $j=1,\dots,n$ and $r \leq 0$;
\item $\Theta^q_{\Psi}(\ltwod) \perp \Theta^q_{\Phi}(\ltwod)$ if $\lambda(E_j^r \cap F_j^r) = 0$ for $j=1,\dots,n$ and $r \leq 0$.
\end{enumerate}
\end{theorem}
\begin{proof}
If $\Theta^q_{\Psi}(\ltwod) \subset \Theta^q_{\Phi}(\ltwod)$, then the range for the analysis operator for the Bessel sequence $\{T_{Xz} \widetilde{D}_{A^r} \psi_j\}$ is contained in that for $\{T_{Yz} \widetilde{D}_{B^r} \phi_j\}$ for $j=1,\dots,n$ and $r \leq 0$.  It follows that the support of
\[ m_j^r(\xi) = \sum_{k \in \Bb{Z}^d} |\hat{\psi}_j(A^{*-r}X^{*-1}(\xi + k))|^2 \]
must be contained in the support of
\[ n_j^r(\xi) = \sum_{k \in \Bb{Z}^d} |\hat{\phi}_j(B^{*-r}Y^{*-1}(\xi + k))|^2 \]
for $j = 1,\dots,n$ and $r \leq 0$.

Items 2 and 3 follow directly.
\end{proof}

\begin{corollary}
In the special case of $X=Y=I$, $A=B$, and $A$ maps the integer lattice into itself, if $m_j^{0}(\xi)$ and $n_j^{0}(\xi)$ have disjoint supports for $j=1,\dots,n$, then $\Cal{U}_{A,I}^q(\Psi)$ and $\Cal{U}_{A,I}^q(\Phi)$ are strongly disjoint.
\end{corollary}
\begin{proof}
If $A$ maps the integer lattice into itself, then for $r \leq 0$, 
\[ \sum_{k \in \Bb{Z}^d} |\hat{\psi}_j(A^{*-r}(\xi + k))|^2 \leq \sum_{k \in \Bb{Z}^d} |\hat{\psi}_j(A^{*-r}\xi + k)|^2. \]
Therefore, for almost every $\xi \in \Bb{R}^d$, if
\[ \sum_{k \in \Bb{Z}^d} |\hat{\psi}_j(A^{*-r}(\xi + k))|^2 = m_j^{r}(\xi) \neq 0, \]
then
\[ \sum_{k \in \Bb{Z}^d} |\hat{\psi}_j(A^{*-r}\xi + k)|^2 = m_j^{0}(A^{*-r}\xi) \neq 0.\]
Hence,
\[ \sum_{k \in \Bb{Z}^d} |\hat{\phi}_j(A^{*-r}\xi + k)|^2 = n_j^{0}(A^{*-r}\xi) =0 \]
and thus
\[ \sum_{k \in \Bb{Z}^d} |\hat{\phi}_j(A^{*-r}(\xi + k))|^2 = n_j^{r}(\xi) = 0. \]
Likewise, if $n_{j}^r(\xi) \neq 0$, then $m_{j}^r(\xi) = 0$.
\end{proof}

\begin{corollary}
If $\Psi$ consists of compactly supported functions and $\Phi$ consists of functions with compactly supported Fourier transforms, then $\Cal{U}^q(\Psi)$ and $\Cal{U}^q(\Phi)$ cannot be similar frames, even if $\Cal{U}(\Psi)$ and $\Cal{U}(\Phi)$ are similar frames.
\end{corollary}
\begin{proof}
If $\Psi$ is compactly supported, then the sums $m_j^r(\xi)$ will supported almost everywhere, whereas if $\Phi$ has compactly supported Fourier transform, then there exists an $r$ with $|r|$ sufficiently large so that $n_j^r(\xi)$ will not be supported almost everywhere.
\end{proof}

\begin{proof}[\bf{PROOF of Theorem 3.}]
The proof of item 1 follows the proof of Theorem \ref{T:W1}.  By analyzing the Bessel sequences $\{T_{Xz}f_j\}$ and $\{T_{Yz}g_j\}$, it follows immediately that $\Theta_{F}(\ltwod) \subset \Theta_{G}(\ltwod)$ only if $\lambda(E_j \setminus F_j) = 0$.  We can likewise analyze the Bessel sequences $\{E_{Al}f_j\}$ and $\{E_{Bl}g_j\}$.  We now have a representation of $\Bb{Z}^d$ on $\overline{span}\{E_{Al}f_j\}$, whose decomposition operator $V:\ltwod \to L^2(\Pi,\lambda,l^2(\Bb{Z}^d))$ is given by
\[ Vf(x)[k] = f(A^{*-1}(x + k)). \]
Immediately we see that $\Theta_{F}(\ltwod) \subset \Theta_{G}(\ltwod)$ only if $\lambda(\tilde{E}_j \setminus \tilde{F}_j) = 0$ for $j=1,\dots,n$.

Item 2 follows from item 1.  Item 3 follows from lemma \ref{L:waveL}.  Finally, for item 4, we compute for any $h_1, h_2 \in \ltwod$:
\begin{multline*}
\sum_{j=1}^{n}\sum_{l \in \Bb{Z}^d}\sum_{z \in \Bb{Z}^d} \langle h_1, E_{Al} T_{Xz} f_j \rangle \overline{\langle h_2, E_{Bl} T_{Yz} g_j \rangle}  \\
= \sum_{j=1}^{n}\sum_{l \in \Bb{Z}^d}\sum_{z \in \Bb{Z}^d} e^{2\pi i \langle Al,Xz \rangle} e^{-2\pi i \langle Bl, Yz\rangle} \langle T_{Xz}^{-1} h_1, E_{Al} f_j \rangle \overline{\langle T_{Yz}^{-1} h_2, E_{Bl} g_j \rangle}.
\end{multline*}
We see that if $\lambda(\tilde{E}_j \cap \tilde{F}_j) = 0$ for $j=1,\dots,n$, then the Bessel sequences $\{E_{Al}f_j\}$ and $\{E_{Bl}g_j\}$ are strongly disjoint, and therefore under the hypothesis $X^{*}A = Y^{*}B$, the Weyl Heisenberg frames are strongly disjoint.
\end{proof}

\subsection{Examples}
We wish to use the above statements to construct examples of strongly disjoint affine frames.

\noindent
{\bf Example 1.} 
Let $\psi$ be a Meyer class wavelet on $\ltwo$.  Recall that the support of $\hat{\psi}$ is contained in the set $[-\frac{4}{3},-\frac{1}{3}) \cup [\frac{1}{3},\frac{4}{3})$ with our version of the Fourier transform.  What we will do is oversample the affine frame at 2 different rates to obtain strongly disjoint Bessel sequences.  Oversampling an affine frame by the factor $N$ is the process of replacing the system
\begin{equation} \label{E:1}
\{\sqrt{2}^n \psi(2^n x - l): n,l \in \Bb{Z}\}
\end{equation}
with the system
\begin{equation} \label{E:2}
\{\frac{\sqrt{2}^n}{\sqrt{N}} \psi(2^n x - \frac{l}{N}): n,l \in \Bb{Z}\}.
\end{equation}
The second oversampling theorem of Chui and Shi states that if (\ref{E:1}) is a tight frame with bound $C_2$ then (\ref{E:2}) is also a tight frame with bound $C_2$ provided $N$ is odd.  We refer to \cite{CS94a, CCMW02a} for detailed information regarding oversampling affine frames.

Thus, we oversample the Meyer wavelet basis by the factors 3 and 13 to obtain two tight affine frames with bound 1.  Moreover, the sums
\[ \sum_{k \in \Bb{Z}} |\hat{\psi}(3(\xi + k))|^2 \text{ and } \sum_{k \in \Bb{Z}} |\hat{\psi}(13(\xi + k))|^2 \]
are supported, when restricted to $[-\frac{1}{2},\frac{1}{2})$, in the sets $[-\frac{4}{9}, -\frac{1}{9}) \cup [\frac{1}{9},\frac{4}{9})$ and $[-\frac{4}{39}, -\frac{1}{39}) \cup [\frac{1}{39},\frac{4}{39})$, respectively.  Hence, the Bessel sequences $\{\frac{1}{\sqrt{3}}\psi(x - \frac{l}{3})\}$ and $\{\frac{1}{\sqrt{13}}\psi(x - \frac{l}{13})\}$ are strongly disjoint, and therefore 
\[ \{\frac{\sqrt{2}^n}{\sqrt{3}} \psi(2^n x - \frac{l}{3})\} \text{ and } \{\frac{\sqrt{2}^n}{\sqrt{13}} \psi(2^n x - \frac{l}{13}) \} \]
are strongly disjoint tight affine frames.

\noindent
{\bf Example 2.}
In example 1, we had to oversample by an odd factor, and so the resulting affine frames were not of the ``standard'' form.  In our second example, we construct strongly disjoint affine frames with the standard form, and moreover, we may have as many strongly disjoint affine frames as desired.

Our example arises from the Frazier-Jawerth frames \cite{FJ90a}, see also \cite{HL00a}.  We construct a frame wavelet whose Fourier transform is supported in the interval $[-\frac{1}{2},-\frac{1}{8}) \cup [\frac{1}{8},\frac{1}{2})$, which satisfies the conditions
\[ \sum_{j \in \Bb{Z}} |\hat{\psi}_0(2^j \xi)|^2 = 1, \qquad \sum_{j \geq 0}\hat{\psi}_0(2^j\xi)\overline{\hat{\psi}_0(2^j(\xi + q))} = 0, \ q \text{ odd,}\]
and which is $C^{\infty}$ (see \cite[Section 5.4]{HL00a} for an example).  It follows that $\psi_0$ is a normalized tight frame wavelet \cite{RS97b}.  Moreover, the function $\psi_n$ defined by $\hat{\psi}_n(\xi) = \hat{\psi}_0(4^n \xi)$ also is a normalized tight frame wavelet.  Finally, by theorem \ref{T:W1}, the affine frames generated by $\psi_n$ are strongly disjoint.

We end, as a result of the proceeding example, with this proposition.
\begin{proposition}
For any integer $N$, there exists a set $\{\psi_1,\dots,\psi_N\}$ of normalized tight frame wavelets which are strongly disjoint such that the Fourier transforms are smooth and compactly supported.
\end{proposition}

\section{acknowledge}
The first and second and fourth authors were supported in part by NSF grants.
The third author's research was supported in part by the
Academic Research Fund No. R-146-000-025-112,
National University of Singapore.
We wish to thank the Mathematics Department at Vanderbilt University for their hospitality during the last three authors' visits there.

\bibliographystyle{amsplain}
\bibliography{biblio}

\end{document}